\input amstex
\magnification=\magstep 2
\documentstyle{amsppt}
\pagewidth{16.9 true cm}
\pageheight{24 true cm}
\TagsOnRight
\def\conv{\mathop{\fam0 conv}}
\def\ess{\mathop{\fam0 ess}}
\topmatter
\title The Kreps-Yan theorem for $L^\infty$
\endtitle
\rightheadtext{The Kreps-Yan theorem}
\author    D.~B.~Rokhlin  \endauthor
\address Zorge str. 5, Department of Mechanics and
         Mathematics, Rostov State University, Rostov-on-Don,
         Russian Federation, 344090
\endaddress
\email  rokhlin\@math.rsu.ru \endemail
\abstract
We prove the following version of the Kreps-Yan theorem.
For any norm closed convex cone $C\subset L^\infty$ such that
$C\cap L_+^\infty=\{0\}$ and $C\supset -L_+^\infty$,
there exists a strictly positive continuous linear functional,
whose restriction on $C$ is non-positive. The proof uses
some tools from convex analysis in contrast to the case of
a weakly Lindel\"of Banach space, where such approach is not needed.
\endabstract
\endtopmatter
\document
\subhead 1. The Kreps-Yan theorem
\endsubhead
Let $\langle X,Y\rangle $ be a pair of Banach spaces in
separating duality [15]. Suppose that
$X$ is endowed with a locally convex topology $\tau$, which is always
assumed to be compatible with this duality, and $K\subset X$ is a
$\tau$-closed convex cone satisfying the condition $K\cap(-K)=\{0\}$.
An element $\xi\in Y$ is called strictly positive
if $\langle x,\xi\rangle>0$ for all $x\in K\backslash\{0\}$.
An element $\xi$ is called non-negaitive if $\langle x,\xi\rangle\ge 0$
for all $x\in K$. We only consider cones $K$ such that the
set of strictly positive functionals is non-empty.

Following [7], we say that the Kreps-Yan theorem
is valid for the ordered space $(X,K)$ with the topology $\tau$
if for any $\tau$-closed convex cone $C$, containing $-K$, the
condition $C\cap K=\{0\} $ implies the existence of a strictly positive
element $\xi\in Y$ such that its restriction on $C$ is non-positive:
$\langle x,\xi\rangle\le 0$, $x\in C$.

If this statement is true for any cone $K$, satisfying the above
conditions, we say that the Kreps-Yan theorem
is valid for the space $(X,\tau)$. It should be mentioned that
in this terminology the Kreps-Yan theorem may be valid for $(X,\tau)$
even if there exists a $\tau$-closed convex cone
$K'\cap(-K')=\{0\}$ such that the set of strictly positive functionals
is empty.

Recall that a space $(X,\tau')$ is said to be Lindel\"of, or
have the Lindel\"of property, if every open cover of $X$ has a
countable subcover [9]. As usual, we denote the weak topology
by $\sigma(X,Y)$.

The next theorem is a modification of theorem 3.1 from
[7], where the more general framework of locally convex
spaces is considered.
\proclaim{Theorem 1} Let $(X,\sigma(X,Y))$ be a Lindel\"of space.
Then Kreps-Yan theorem is valid for the space $(X,\tau)$.
\endproclaim
\demo{Proof} Let $x\in K\backslash\{0\}$, then
$x\not\in C$ and by the separation theorem [15, th. II, 9.2])
there exists an element $\xi_x\in Y$ such that
$$ \langle y,\xi_x\rangle<\langle x,\xi_x\rangle,\ \ y\in C.$$
But $C$ is a cone, hence we get the inequality
$\langle y,\xi_x\rangle\le 0$, $y\in C$. In addition,
$-K\subset C$. Consequently,
$$ \langle x,\xi_x\rangle>0,\ \ \
   \langle z,\xi_x\rangle\ge 0,\ \ z\in K.$$

Consider the following family of sets
$$ A_x=\{y\in X:\langle y,\xi_x\rangle>0\},\ \ x\in K\backslash\{0\}$$
and let $A_0=\{y\in X:|\langle y,\eta\rangle|<1\}$, where
$\eta$ is a strictly positive functional.
The sets $A_x$, $x\in K$ are open in the topology $\sigma(X,Y)$
and constitute an open cover of $K$.  Moreover, the cone $K$
is closed in $\sigma(X,Y)$, because all topologies compatible
with the duality $\langle X,Y\rangle$ have the identical collection
of closed convex sets. In view of Lindel\"of property, this implies
the existence of the following countable subcover:
$K\subset\cup_{i=0}^\infty A_{x_i}$, where $x_0=0$.

Let $\alpha_i=1/(\|\xi_{x_i}\| 2^i)$, then
$\sum_{i=1}^\infty\alpha_i \xi_{x_i} $ converges in the norm topology
to some element $\xi\in Y$. Evidently, $\xi\le 0$ on $C$.
Moreover, $\xi$ is strictly positive. Indeed,
for any element $x\in K\backslash\{0\}$ there exists
a $\lambda>0$ such that $\lambda x\neq A_0$.
Consequently, $\lambda x\in A_{x_k}$ for some $k\ge 1$ and
$$ \langle\lambda x,\xi\rangle=\sum_{i=1}^\infty
   \alpha_i\langle \lambda x,\xi_{x_i}\rangle\ge
   \alpha_k\langle\lambda x,\xi_{x_k}\rangle>0. $$
This completes the proof.
\enddemo

A similar result was proved in [7] under the next
condition, conceptually connected with the Halmos-Savage theorem
[8, lem. 7]. For any family of non-negative functionals
$\{\xi\}_{\beta\in I}\subset Y$ there exist a countable subset
$\{\xi_{\beta_i}\}_{i=1}^\infty$ with the following property: if
for $x\in K\backslash\{0\}$ there exists a $\beta\in I$ such that
$\langle x,\xi_\beta\rangle>0$, then
$\langle x,\xi_{\beta_i}\rangle>0$ for some $i$.

We prefer to require that the space $(X,\sigma(X,Y))$ verifies
the more standard Lindel\"of condition. Evidently, this condition
is satisfied if any topology of the space $X$, compatible with
the duality $\langle X,Y\rangle$, has the Lindel\"of property.

Also note, that the space $X$ is Lindel\"of if it may be
represented as the union of a countable collection of compact sets.
Hence, a reflecsive space $X$ is Lindel\"of in the weak topology
$\sigma(X,X^*)$ (shortly: weakly Lindel\"of)  in view of the weak
compactness of the unit ball, and the space $X^*$ is Lindel\"of in
the $*$@-weak topology $\sigma(X^*,X)$ by the Banach-Alaoglu theorem.
Consequently, the Kreps-Yan theorem is valid for any reflecsive space
with the norm topology and for the space $(X^*,\sigma(X^*,X))$.

A Banach space $X$ is called weakly compactly generated
(shortly: WCG), if $X$ contains a weakly compact subset whose
linear span is dense in $X$. Corson conjectured that
the notions of weakly Lindel\"of and WCG spaces are equivalent [1].
The one-half of this conjecture was confirmed in [16] (see also
[5, th. 12.35]): every WCG space is weakly Lindel\"of (the converse
implication appeared to be false in general as follows from [13], [11])).
Therefore, the Kreps-Yan theorem is valid for any WCG space,
endowed with the norm topology.

Let $(\Omega,\Cal F,\bold P)$ be a probability space. Denote by
$L^p=L^p(\Omega,\Cal F,\bold P)$, $1\le p\le\infty$
the Banach spaces of equivalence classes of measurable functions,
whose $p$'s power is integrable if $1\le p<\infty$ (respectively,
which are essentially bounded if $p=\infty$). The above arguments
imply the following result (compare with [7], [14, th. 1.4]): the
Kreps-Yan theorem is valid for the spaces $(L^p,\tau_p)$,
$1\le p\le\infty$, where $\tau_p$ is the norm topology for
$1\le p<\infty$, and $\tau_\infty$ is the weak-star topology
$\sigma(L^\infty,L^1)$.

Indeed, the spaces $L^p$, $1<p<\infty$ are reflecsive, the
topology $\sigma(L^\infty,L^1)$ of the space $L^\infty$ coincides
with the weak-star one, and $L^1$ is a WCG space [3, p.143].

On the other hand, it is shown in [7, example 2.1] that the Kreps-Yan
theorem may fail even if $(X,K)$ is a Banach lattice
(with the norm topology).
So, the imposed Lindel\"of condition is not superfluous.
Note also, that theorem 1 does not imply the
validity of the Kreps-Yan theorem for the space $L^\infty$ with
the norm topology: it is known that even the space of bounded
sequences is not weakly Lindel\"of [1, example 4.1(i)].

\subhead 2. The case of $L^\infty$
\endsubhead
Let $L_+^\infty$ be the cone, generating the natural
order structure on $L^\infty$. Our main result is the following.
\proclaim{Theorem 2} The Kreps-Yan theorem is valid for the
ordered space $(L^\infty,L_+^\infty)$ with the norm topology.
\endproclaim

Let $C\subset L^\infty$ be a norm closed convex cone,
satisfying the conditions
$$C\cap L^\infty_+=\{0\},\ \ -L^\infty_+\subset C.$$
Put
$C_\varepsilon=\{x\in C:\ess\inf x\ge-\varepsilon\}.$

Recall that the dual of $L^\infty$ (with the norm topology)
coincides with the Banach
space $ba=ba(\Omega,\Cal F,\bold P)$ of all bounded finitely
additive measures $\mu$ on $(\Omega,\Cal F)$ with the property
that $\bold P (A)=0$ implies $\mu(A)=0$ ([4]).
Let
$$  ba^+=\{\mu\in ba:\langle x,\mu\rangle\ge 0,\ x\in L_+^\infty\}$$
be the set of non-negative elements of $ba$. A probability
measure $\bold Q$ is identified with the continuous functional on
$L^\infty$ by the formula
$ \langle x,\bold Q\rangle=\int_\Omega x\,d\bold Q.$

Denote by $\conv A$ the convex hull of the set $A$.
If $D\subset \Omega$ we put
$$ I_D(\omega)=1,\ \omega\in D;\ \ I_D(\omega)=0,\ \omega\not\in D.$$
The expectation with respect to $\bold P$ is denoted by $\bold E x$.

\proclaim{Lemma} There exists a probability measure
$\bold Q$ equivalent to $\bold P$ such that
$$   \sup_{x\in C_1} \langle x,\bold Q\rangle<\infty.$$
\endproclaim
\demo{Proof} It suffice to show that the set $C_1$ is
bounded in probability, since this reduces the assertion
of lemma to the Yan theorem [17] (see also [12, p.145]).
We literally follow the argumentation of [2, prop.~3.1],
where somewhat special analog of the set $C$ is considered.

Let the set $C_1$ be unbounded in probability. Then there exist
a sequence of elements $x_n\in C_1$, $n\ge 1$ and a number $\alpha>0$
such that $\bold P(x_n\ge n)>\alpha$. Elements
$y_n=\min\{x_n/n,1\}$ belong to $C_{1/n}\subset C_1$ and
$$ \bold P(y_n=1)=\bold P(x_n/n\ge 1\}>\alpha.$$
By [2, lem. A1.1] there exist a sequence
$$z_n\in\conv(y_n,y_{n+1},\dots)\subset C_{1/n},$$
converging a.s. to $z:\Omega\mapsto [0,1]$.
Further, the inequality
$$\bold E y_n \ge \bold E(y_n I_{y_n})\ge
\bold E(I_{\{y_n=1\}})-\bold E(I_{\{y_n<1\}}/n)\ge
  \alpha-1/n, $$
implies that $\bold E z_n\ge\alpha-1/n$ and by Lebesgue's
dominated convergence theorem
$$ \bold E z=\lim_{n\to\infty}\bold E z_n\ge \alpha.$$
Hence
$$\bold P(z>0)=\beta\ge\bold E(zI_{\{z>0\}})=\bold E z\ge\alpha.$$
By Egorov's theorem $z_n\to z$ uniformly on a set
$\Omega'$: $\bold P(\Omega')\ge 1-\beta/2$.
The functions $w_n=\min\{z_n,I_{\Omega'}\}$ belong to $C$ and
$ w_n=z_n I_{\Omega'}\to z I_{\Omega'}$ in the norm topology
of $L^\infty$. We obtain a contradiction, since
$$ \bold P(z I_{\Omega'}>0)=\bold P(\Omega')+\bold P(z>0)-
   \bold P(\Omega'\cup \{z>0\})\ge\beta/2,$$
This completes the proof.
\enddemo

Now we need some additional notation, used in convex analysis
(e.g. [10]). Let again $\langle X,Y\rangle$ be
a pair of Banach spaces in duality. The indicator and support
functions of a convex set $A\subset X$ are defined by the formulas
$$ \delta A (x)=0,\ x\in A,\ \ \delta A(x)=+\infty,\ x\not\in A;\ \ \
   sA(\xi)=\sup_{x\in A}\langle x,\xi\rangle. $$
The same notation is used if $A\subset Y$. The sets
$$ A^\circ=\{\xi\in Y:\langle x,\xi\rangle\le 1,\ \ x\in A\},\ \ \
   A^{\circ\circ}=\{x\in X:\langle x,\xi\rangle\le 1,\ \ \xi\in A^\circ\}$$
are called polar and bipolar of $A$.

The  Young-Fenchel transform of a function
$f:X\mapsto [-\infty,+\infty]$ is defined as follows
$$ f^*(\xi)=\sup_{x\in X}(\langle x,\xi\rangle-f(x))$$
The function
$$ (f_1\oplus f_2)(x)=\inf\{f(x_1)+f(x_2):x_1+x_2=x\}$$
is called an infimal convolution of $f_1$, $f_2$.

Note, that the support function of a set $A$ is equal to the
Minkowski function $\mu A^\circ$ of the polar $A^\circ$:
$$ sA(\xi)=\mu A^\circ(\xi),\ \
   \mu A^\circ=\inf\{\lambda>0: x\in\lambda A^\circ\}.$$
We will use the next formula for its Young-Fenchel transform:
$$ \split
   (\mu A^\circ)^*(x) &= \sup_{\xi\in Y}(\langle x,\xi\rangle-
   \inf\{\lambda>0: \xi\in\lambda A^\circ\})=
   \sup_{\lambda>0}\sup_{\xi\in\lambda A^\circ}
   (\langle x,\xi\rangle-\lambda)\\
   &=\sup_{\lambda>0}\lambda(\sup_{\eta\in A^\circ}
    \langle x,\eta\rangle-1)=\delta A^{\circ\circ}(x)
   \endsplit$$

\demo{Proof of the theorem 2} Let $\bold Q$ be a measure, introduced
in the previous lemma. Put
$$ \varphi(\varepsilon)=-\sup_{x\in C_{\varepsilon}}
    \langle x,\bold Q\rangle. $$
Note, that $C_\varepsilon=\emptyset$ for $\varepsilon<0$,
$C_0=0$, and $C_\varepsilon=\varepsilon C_1$ for $\varepsilon>0$.
Since the support function of an empty set is equal to $-\infty$,
we get
$$ \varphi(\varepsilon)=-\varepsilon\varphi(1)
   +\delta(-\infty,0)(\varepsilon), \ \ \varphi(1)\le 0.$$

Denote by $\Cal P$ the set of all probability measures $\bold P'$,
absolutely continuous with respect to $\bold P$. We have
$$ \ess\inf x=\inf_{\bold P'\in \Cal P}\langle x,\bold P'\rangle
   =-s(-\Cal P)(x)$$
and $C_\varepsilon=C\cap\{x\in L^\infty:s(-\Cal P)(x)\le\varepsilon\}$.
So, for $\tau<0$ the function $\varphi^*$ has the following
representation:
$$ \split
   \varphi^*(\tau) &=\sup_{\varepsilon\ge 0}\sup_{x\in C_{\varepsilon}}
   (\varepsilon\tau+\langle x,\bold Q\rangle)=\sup_{x\in C}
   \sup_{\varepsilon\ge 0}\sup_{s(-\Cal P)(x)\le\varepsilon}
   (\varepsilon\tau+\langle x,\bold Q\rangle)\\
   &=\sup_{x\in C}(\tau\cdot s(-\Cal P)(x)+\langle x,\bold Q\rangle).
   \endsplit$$
For $\lambda=-\tau$ we obtain
$$ \split
   \varphi^*(-\lambda) &=\sup_{x\in L^\infty}(\langle x,\bold Q\rangle
   -\lambda\cdot s(-\Cal P)(x)-\delta C(x))=
    (s(-\lambda \Cal P)+\delta C)^*(\bold Q)\\
     &=((s(-\lambda \Cal P))^*\oplus (\delta C)^*)(\bold Q).
    \endsplit$$
The last equality uses the formula
$ (f_1+f_2)^*=f_1^*\oplus f_2^*$
which is valid, in particular, if one of the functions is continuous
on the whole space [6]. In our case it is true for
$s(-\lambda \Cal P)$.

Using the identities
$$ s(-\lambda \Cal P)^*=(\mu(-\lambda \Cal P)^\circ)^*=
\delta(-\lambda \Cal P)^{\circ\circ},
  \ \ (\delta C)^*=sC=\delta C^\circ,$$
we get
$$ \varphi^*(-\lambda)=
    (\delta(-\lambda \Cal P)^{\circ\circ}\oplus \delta C^\circ)(\bold Q)
    =\delta((-\lambda \Cal P)^{\circ\circ}+C^\circ)(\bold Q).$$

On the other hand, the Young-Fenchel transform of $\varphi$
is given by
$$ \varphi^*(\tau)=
  \sup_\varepsilon(\varepsilon\tau-\varphi(\varepsilon))
  =\sup_{\varepsilon\ge 0}(\varepsilon(\tau-\varphi(1))=
   \delta(-\infty,\varphi(1)](\tau).$$
It follows that $\varphi^*(-\lambda)=0$
for $\lambda>-\varphi(1)$. Thus,
$$ \bold Q\in C^\circ+(-\lambda \Cal P)^{\circ\circ},\ \
   \lambda\in (-\varphi(1),+\infty)$$
and there exists an element $\mu\in C^\circ$ such that
$$ \mu=\bold Q+\nu,\ \ \nu\in -(-\lambda \Cal P)^{\circ\circ}
                       =\lambda \Cal P^{\circ\circ}. $$

But $\Cal P^{\circ\circ}$ coincides with the $\sigma(ba,L^\infty)$ closed
convex hull of the set $\Cal P\cup\{0\}\subset ba^+$ [15, th. IV, 1.5].
Hence, $\nu\in ba^+$ and $\mu$ is a desired functional: it is
strictly positive and $\langle x,\mu\rangle\le 0$, $x\in C$.
The proof is complete.
\enddemo

Finally, we mention that the case of $L^\infty$ with the norm
topology is of special interest for mathematical finance in
view of characterization of the No Free Lunch with Vanishing Risk
condition [2]. We hope to discuss such applications elsewhere.

Another interesting question concerns the precise description
of Banach lattices such that the Kreps-Yan theorem is valid
under the norm topology.

\enddocument